\documentclass[11pt]{article}
\usepackage{amsmath,amssymb,amsthm,bm,graphicx,color,epsfig}

\newtheorem{theorem}{Theorem}[section]
\newtheorem{lemma}[theorem]{Lemma}

\newcommand{\R}{\mathbb{R}}

\renewcommand{\i}{\imath}

\newcommand{\eps}{\varepsilon}
\newcommand{\p}{\partial}

\begin{document}

\title{Sparse Fourier Transform via Butterfly Algorithm}
\author{Lexing Ying\\
  Department of Mathematics, University of Texas, Austin, TX 78712}

\date{December 2007}
\maketitle

\begin{abstract}
  We introduce a fast algorithm for computing sparse Fourier
  transforms supported on smooth curves or surfaces. This problem
  appear naturally in several important problems in wave scattering
  and reflection seismology. The main observation is that the
  interaction between a frequency region and a spatial region is
  approximately low rank if the product of their radii are bounded by
  the maximum frequency. Based on this property, equivalent sources
  located at Cartesian grids are used to speed up the computation of
  the interaction between these two regions. The overall structure of
  our algorithm follows the recently-introduced butterfly algorithm.
  The computation is further accelerated by exploiting the
  tensor-product property of the Fourier kernel in two and three
  dimensions. The proposed algorithm is accurate and has an $O(N \log
  N)$ complexity.  Finally, we present numerical results in both two
  and three dimensions.
\end{abstract}




\section{Introduction}
\label{sec:intro}

We consider the rapid computation of the following sparse Fourier
transform problem. Let $N$ be a large integer and $X$ and $K$ be two
smooth (or piecewise smooth) curves in the unit box $[0,1]^2$.
Suppose $\{ x_i, i\in I\}$ and $\{k_j, j\in J\}$ are, respectively,
the samples of $NX$ and $NK$, where $NX := \{N\cdot p, p \in X\}$ and
$NK$ is defined similarly. Given the sources $\{f_j, j\in J\}$ at
$\{k_j, j \in J\}$, the problem is to compute the potentials $\{u_i, i
\in I \}$ defined by
\begin{equation}
  u_i = \sum_{j\in J} e^{ 2\pi\i x_i \cdot k_j / N } f_j,
  \label{eq:sfft}
\end{equation}
where $\i = \sqrt{-1}$. In most cases, $\{x_i, i\in I\}$ and $\{k_j,
j\in J\}$ sample $NX$ and $NK$ with a constant number of samples per
unit length. As a result, $\{x_i, i\in I\}$ and $\{k_j, j\in J\}$ are
of size $O(N)$.  A similar problem can be defined in three dimensional
space where $X$ and $K$ are smooth surfaces in $[0,1]^3$. However, our
discussion here focuses on the two dimensional case, as the algorithm
can be copied verbatim to the three dimensional case.

Direct evaluation of \eqref{eq:sfft} clearly requires $O(N^2)$ steps,
which can be quite expensive for large values of $N$. In this paper,
we propose an $O(N \log N)$ approach based on the butterfly algorithm
\cite{michielssen-1996-mmda,oneil-2007-ncabft}. Our algorithm starts
by generating two quadtrees $T_X$ and $T_K$ for $\{x_i\}$ and
$\{k_j\}$, respectively, where each of their leaf boxes is of unit
size.  The main observation is the following low rank property. Let
$A$ and $B$ be are two boxes from $T_X$ and $T_K$, respectively. If
the widths $w^A$ and $w^B$ of these two boxes satisfy the condition
$w^A \cdot w^B \le N$, then the interaction $e^{2\pi \i x\cdot k/N}$
between $x\in A$ and $k\in B$ is approximately of low rank.  This
property implies that one can reproduce the potential in $A$ with a
set of equivalent charges whose degree of freedom is bounded by a
constant independent of $N$.

The algorithm first constructs the equivalent charges for the leaf
boxes of $T_K$. Next, we traverse up $T_K$ to construct the equivalent
charges of the non-leaf boxes from the ones of their children. For
each non-leaf box $B$, we construct a set of equivalent charges for
each box $A$ in $T_X$ that satisfies $w^A = N/w^B$. At the end of this
step, we hold the equivalent charges of the root box of $T_K$ for each
leaf box of $T_X$. Finally, we visit all of the leaf boxes of $T_X$
and utilize the equivalent charges of the root box of $T_K$ to compute
the potentials $\{u_i, i\in I \}$.

The sparse Fourier transforms in \eqref{eq:sfft} appears naturally in
several contexts. One example is from the calculation of the
time-harmonic scattering field \cite{colton-1983-iemst}. Suppose that
$D$ is a smooth surface and $N$ is the wave number. The scattering
field $u$ satisfies the Helmholtz equation $-N^2 u - \Delta u = 0$ in
$\R^d \setminus D$. Its far field pattern $u^\infty(\hat{x})$ for
$\hat{x}$ on the unit sphere, which is highly important for many
scattering problems, is defined by
\begin{equation}
  u_\infty (\hat{x}) = \int_{\p D} e^{ -\i N \hat{x}\cdot y } f(y) d s(y),
  \label{eq:far}
\end{equation}
where $f$ is some function supported on the boundary $\p D$. After we
rescale $\hat{x}$ and $y$ by a factor of $N$, \eqref{eq:far} takes the
form of \eqref{eq:sfft}.  Another example of \eqref{eq:sfft} appears
in the depth stepping algorithm in reflection seismology
\cite{fomel-2007-fcpft}.

The rest of this paper is organized as follows. In Section
\ref{sec:algo}, we prove the main analytic result and describe our
algorithm in detail. After that, we provide numerical results for both
the two and three dimensional cases in Section \ref{sec:num}.
Finally, we conclude in Section \ref{sec:conc} with some discussions
on future work.

\section{Algorithm Description}
\label{sec:algo}

The main theoretical component of our algorithm is the following
theorem. Following the discussion in Section \ref{sec:intro}, we use
$T_X$ and $T_K$ to denote the quadtrees generated from $X$ and $K$,
respectively.

\begin{theorem}
  \label{thm:main}
  Let $A$ be a box in $T_X$ and $B$ be a box in $T_K$. Suppose the
  width $w^A$ and $w^B$ of $A$ and $B$ satisfy $w^A w^B = N$. Then,
  for any $\eps > 0$, there exists a constant $T(\eps)$ and functions
  $\{\alpha_t(x), 1\le t \le T(\eps) \}$ and $\{\beta_t(k), 1\le t\le
  T(\eps)\}$ such that
  \[
  \left| 
    e^{2\pi\i x\cdot k/N} - \sum_{t=1}^{T(\eps)} \alpha_t(x) \beta_t(k)
  \right|
  \le \eps
  \]
  for any $x \in A$ and $k \in B$.
\end{theorem}

The proof of Theorem \ref{thm:main} is based on the following
elementary lemma.
\begin{lemma}
  For any $Z>0$ and $\eps>0$, let $S = \lceil \max(4e\pi Z, \log_2
  (1/\eps)) \rceil$. Then
  \[
  \left|    e^{2\pi\i x} - \sum_{t=0}^{S-1} \frac{(2\pi\i x)^t}{t!}  \right|  \le \eps
  \]
  for any $x$ with $|x| \le Z$.
\end{lemma}

\begin{proof}[Proof of Theorem \ref{thm:main}]
  Let us use $c^A = (c^A_1, c^A_2)$ and $c^B = (c^B_1, c^B_2)$ to
  denote the lower left corners of boxes $A$ and $B$, respectively.
  Writing $x=c^A+x'$ and $k=c^B+k'$, we have
  \[
    e^{2\pi\i x\cdot k/N} = e^{2\pi\i(c^A+x')\cdot(c^B+k')/N} 
    = 
    e^{2\pi\i (c^A+x')\cdot c^B /N} \cdot 
    e^{2\pi i x'\cdot k'/N} \cdot 
    e^{2\pi\i c^A \cdot k'/N}.
  \]
  Notice that the first and the third terms depends only on $x'$ and
  $k'$, respectively. Therefore, we only need to construct a
  factorization for the second term. Since $|x'|\le \sqrt{2}w^A$ and
  $|k'|\le \sqrt{2} w^B$, $|x'\cdot k'| /N \le 2$. Invoking the lemma
  for $Z=2$, we obtain the following approximation with $S(\eps)$
  terms:
  \[
  \left| e^{2\pi\i x'\cdot k'/N} - \sum_{t=0}^{S(\eps)-1} \frac{(2\pi\i x'\cdot k'/N)^t}{t!}  \right|  \le \eps.
  \]
  After expanding each term of the sum using $x'\cdot k' = (x'_1 k'_1
  + x'_2 k'_2)$, we have an approximate expansion 
  \[
  \left|  e^{2\pi\i x\cdot k/N} - \sum_{t=1}^{T(\eps)} \alpha_t(x) \beta_t(k)  \right|  \le \eps
  \]
  where $T(\eps)$ only depends on the accuracy $\eps$.
\end{proof}

\subsection{Equivalent sources}
\label{subsec:algo_eqn}
Given two boxes $A$ and $B$ with $w^A w^B = N$, we denote $u^{AB}(x)$
the potential field in $A$ generated by the charges inside $B$:
\[
u^{AB}(x) = \sum_{j: k_j \in B} e^{ 2\pi\i x_i \cdot k_j / N } f_j.
\]
The theorem implies that the field $u^{AB}$ can be approximately
reproduced by a group of carefully selected {\em equivalent charges}
inside $B$. For the efficiency reason to be discussed shortly, we pick
these charges to be located on a Cartesian grid in $B$,
\[
\left\{ k^B_{lm} := \left( c^B_1 + l \frac{w^B}{p-1}, c^B_2 + m \frac{w^B}{p-1} \right), \quad l,m=0,1,\cdots p-1 \right\},
\]
where $p$ is a constant whose value depends on the prescribed accuracy
$\eps$. The corresponding equivalent charges are denoted by $\{
f^{AB}_{lm} \}$. To construct $\{ f^{AB}_{lm} \}$, we first select a
group of points
\[
\left\{ x^A_{lm} := \left( c^A_1 + l \frac{w^A}{p-1}, c^A_2 + m \frac{w^A}{p-1} \right), \quad l,m=0,1,\cdots p-1 \right\}.
\]
located on a Cartesian grid in $A$. $\{ f^{AB}_{lm} \}$ are computed
by ensuring that they reproduce the they generate the field
$\{u^{AB}_{lm} := u^{AB}(x^A_{lm})\}$ at the points $\{ x^A_{lm}\}$;
i.e.,
\[
\sum_{l'm'} e^{2\pi\i x^A_{lm} \cdot k^B_{l'm'}/N} f^{AB}_{l'm'} =  u^{AB}_{lm}.
\]
Writing this into a matrix form $M f = u$ and using the definitions of
$\{x^A_{lm}\}$ and $\{ k^B_{l'm'}\}$, we can decompose the $p^2 \times
p^2$ matrix $M$ into a Kronecker product $M = M_1 \otimes M_2$, where
\[
M_1 = 
\left(
  e^{2\pi\i \left(c_1^A + l\frac{w^A}{p-1}\right) \left(c_1^B + l'\frac{w^B}{p-1}\right)/N } 
\right)_{ll'}\quad
M_2 = \left(
  e^{2\pi\i \left(c_2^A + m\frac{w^A}{p-1}\right) \left(c_2^B + m'\frac{w^B}{p-1}\right)/N }
\right)_{mm'}.
\]
Since $(M_1 \otimes M_2)^{-1} = M_1^{-1} \otimes M_2^{-1}$, in order
to compute $f = M^{-1} u$ we only need to invert the $p \times p$
matrices $M_1$ and $M_2$. Expanding the formula for $M_1$, we get
\begin{equation}
(M_1)_{ll'} = e^{2\pi\i \left(c^A_1 + l\frac{w^A}{p-1}\right) c^B_1 /N} \cdot
e^{2\pi\i \frac{l l'}{(p-1)^2} \frac{w^A w^B}{N} } \cdot
e^{2\pi\i c^A_1 \left(l'\frac{w^B}{p-1}\right) / N}
\label{eq:M1}
\end{equation}
Noticing that the first and the third terms depend only on $l$ and
$l'$, respectively, and $w_A w_B/N = 1$, we can rewrite $M_1$ into a
factorization $M_1 = M_{11} \cdot G \cdot M_{12}$, where $M_{11}$ and
$M_{12}$ are two diagonal matrices and the center matrix $G$ given by
$(G)_{ll'} = e^{2\pi\i \frac{l l'}{(p-1)^2}}$ is independent of $N$
and the boxes $A$ and $B$. The situation for $M_2$ is exactly the
same. The result of this discussion is that we reduce the complexity
of $f = M^{-1} u$ from $O(p^4)$ to $O(p^3)$ using the Kronecker
product structure of $M$. In fact, one can further reduce it to
$O(p^2\log p)$ since the matrix $G$ is a fractional Fourier transform
(see \cite{bailey-1991-ffta}).

\begin{figure}[h]
 \begin{center}
   \includegraphics[height=1.5in]{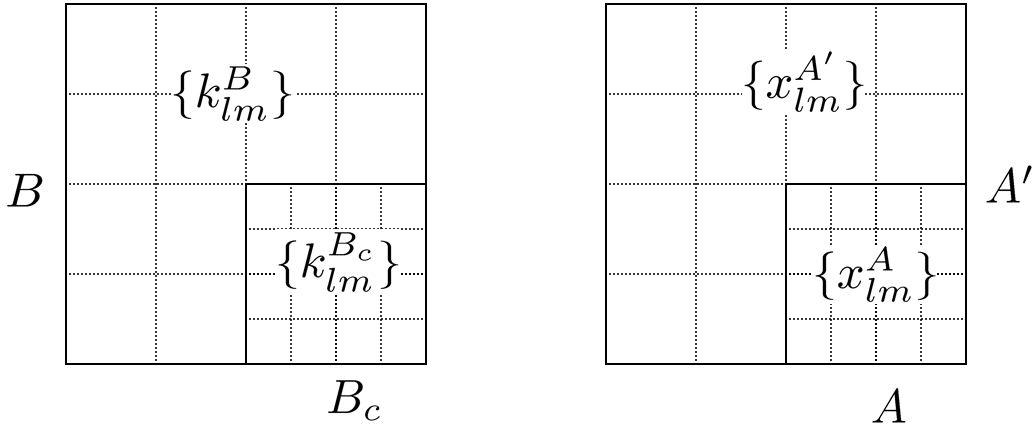}
  \end{center}
  \caption{The construction of $\{f^{AB}_{lm}\}$ using
    $\{f^{A'B_c}_{lm}\}$. $B_c$ is one of $B$'s child boxes and $A'$
    is $A$'s parent box. We first evaluate the potentials
    $\{u^{AB}_{lm}\}$ located at $\{x^A_{lm}\}$ using
    $f^{A'B_c}_{l'm'}$ and then find $\{f^{AB}_{lm}\}$ so that they
    produce the same potentials.  }
  \label{fig:cross}
\end{figure}

The procedure we just described fails to be efficient when we compute
$\{f^{AB}_{lm}\}$ for a large box $B$, as typically $B$ contains a
large number of points $\{k_j\}$ and this makes the evaluation of
$\{u^{AB}_{lm}\}$ quite expensive. The second ingredient of the
butterfly algorithm addresses this problem. In our setting, suppose
that $B$ is a non-leaf box of $T_K$, $A$ is a box in $T_X$, and $w^A
w^B = N$. We denote the children of $B$ by $B_c, c=1,\cdots,4$ and the
parent of $A$ by $A'$. Suppose that one constructs the equivalent
charges in a bottom-up traversal of $T_K$. Hence, when one reaches
$B$, the equivalent charges of $B_c$ have already been computed. The
idea is then to use the equivalent charges $\{f^{A'B_c}_{lm}\}$ of
$B_c$ to compute the check potentials $\{u^{AB}_{lm}\}$ since
$A\subset A'$; i.e.,
\[
u^{AB}_{lm} \approx \sum_{c=1}^4 \left(
  \sum_{l'm'} e^{2\pi \i x^A_{lm}\cdot k^{B_c}_{l'm'}/N} f^{A'B_c}_{l'm'}
\right)
\]
for any $l,m$. The inner sum $\sum_{l'm'} e^{2\pi \i x^A_{lm}\cdot
  k^{B_c}_{l'm'}/N} f^{A'B_c}_{l'm'}$ for each fixed $i$ can be
rewritten into a matrix form $E f$. Using again the Kronecker
product, we can decompose $E$ as $E_1 \otimes E_2$ where
\[
E_1 = 
\left(
  e^{2\pi\i \left(c_1^A + l\frac{w^A}{p-1}\right) \left(c_1^{B_c} + l'\frac{w^{B_c}}{p-1}\right)/N } 
\right)_{ll'}\quad
E_2 = \left(
  e^{2\pi\i \left(c_2^A + m\frac{w^A}{p-1}\right) \left(c_2^{B_c} + m'\frac{w^{B_c}}{p-1}\right)/N }
\right)_{mm'}.
\]
Expanding the formula for $E_1$, we get 
\begin{equation}
(E_1)_{ll'} =  e^{2\pi\i \left(c^A_1 + l\frac{w^A}{p-1}\right) c^{B_c}_1 /N } \cdot
e^{2\pi\i \frac{l l'}{(p-1)^2} \frac{w^A w^{B_c}}{N} } \cdot
e^{2\pi\i c^A_1 \left(l'\frac{w^{B_c}}{p-1}\right) /N }
\label{eq:E1}
\end{equation}
Noticing that the first and the third terms depend only on $l$ and
$l'$ respectively and $w^A w^{B_i}/N = 1/2$, we can write $E_1$ into a
factorization $E_1 = E_{11} \cdot H \cdot E_{12}$ where $E_{11}$ and
$E_{12}$ are again diagonal matrices and the matrix $H$ given by
$(H)_{ll'} = e^{\pi\i \frac{l l'}{(p-1)^2}}$ is independent of $N$.

\subsection{Algorithm}

We now give the overall structure of our algorithm. It contains the
following steps:
\begin{enumerate}
\item Construct the quadtrees $T_X$ and $T_K$ for the point sets $\{x_i,
  i\in I\}$ and $\{k_j, j\in J\}$, respectively. These trees are
  constructed adaptively and all the leaf boxes are of unit size.
\item Construct the equivalent charges for the leaf boxes in $T_K$.
  Suppose that $A_0$ is the root box of $T_X$. For each leaf box $B$
  in $T_K$, we calculate $\{f^{A_0 B}_{lm}\}$ by matching the
  potentials $\{u^{A_0 B}_{lm}\}$.
\item Travel up in $T_K$ and construct the equivalent charges for the
  non-leaf boxes in $T_K$. For each non-leaf box $B$ in $T_K$ and each
  box $A$ in $T_X$ with width $w^A = N/ w^B$, we construct
  $\{f^{AB}_{lm}\}$ from the equivalent charges of
  $\{f^{A'B_c}_{lm}\}$ where $A'$ is $A$'s parent and
  $\{B_c,c=1,\cdots,4\}$ are the $B$'s children.
\item Compute $\{u_j, j\in J\}$. Let $B_0$ be the root box of $T_K$.
  For each leaf box $A$ in $T_X$ and each $j$ such that $x_j \in A$,
  we approximate $u_j$ with
  \[
  \sum_{lm} e^{2\pi \i x_j\cdot k^{B_0}_{lm}/N} f^{A B_0}_{lm}.
  \]
\end{enumerate}

Let us first estimate the number of operations of the proposed
algorithm. The first step clearly takes only $O(N \log N)$ operations
since there are at most $O(N)$ points in both $\{x_i, i\in I\}$ and
$\{k_j, j\in J\}$. By assumption, both $X$ and $K$ are smooth (or
piecewise smooth) curves in $[0,1]^2$. Therefore, for a given size
$w$, there are at most $O(N/w)$ non-empty boxes in both $T_X$ and
$T_K$. In particular, we have at most $O(N)$ leaf boxes in $T_X$ and
$T_K$. This implies that the second and fourth steps take at most
$O(N)$ operations. To analyze the third step, we estimate level by
level. For a fixed level $t$, there are at most $2^t$ non-empty boxes
in $T_K$ on that level, each of size $N/2^t$. For each box $B$ on
level $t$, we need to construct $\{f^{AB}_{lm}\}$ for all the boxes
$A$ in $T_X$ with size $N/(N/2^t) = 2^t$. It is clear that there are
at most $N/2^t$ non-empty boxes in $T_X$ of that size. Since the
construction for each set of equivalent charges take only constant
operations, the total complexity for level $t$ is $O(2^t \times N/2^t)
= O(N)$. As there are at most $O(\log N)$ levels, the third step takes
at most $O(N \log N)$ operations.  Summing over all the steps, we
conclude that our algorithm is $O(N \log N)$.

Our algorithm is also efficient in terms of storage space. Since we
give explicit construct formulas \eqref{eq:M1} and \eqref{eq:E1} for
constructing the equivalent charges, we only need to store the
equivalent charges during the computation. This is where our algorithm
differs from the one in \cite{oneil-2007-ncabft} where they need to
require $O(N\log N)$ small matrices for the interpolative
decomposition \cite{cheng-2005-clrm,gu-1996-eacsrqf}. As we mentioned
early, at each level, we need to keep $O(N)$ equivalent charges, each
of which takes $O(1)$ storage space. Noticing that, at any point of
the algorithm, we only need to keep the equivalent charges for two
adjacent levels, therefore the storage requirement of our algorithm is
only $O(N)$.

The three dimensional case is similar. Since $X$ and $K$ are smooth
surfaces in $[0,1]^3$, the number of points in $\{x_i, i\in I\}$ and
$\{k_j, j\in J\}$ are $O(N^2)$ instead. The Kronecker product
decomposition is still valid and, therefore, we can construct the
equivalent charges efficiently in $O(p^4)$ (or even $O(p^3 \log p)$)
operations instead of $O(p^6)$. The algorithm remains exactly the same
and a similar complexity analysis gives an operation count of
$O(N^2\log N)$, which is almost linear in terms of the number of
points.

\section{Numerical Results}
\label{sec:num}

In this section, we provide some numerical examples to illustrate the
properties of our algorithm. All of the numerical results are obtained
on a desktop computer with a 2.8GHz CPU. The accuracy of the algorithm
depends on $p$, which is the size of the Cartesian grid used for the
equivalent charges. In the following examples, we pick $p=5,7$, or
$9$. The larger the value of $p$, the better the accuracy.

\subsection{2D case}

For each example, we sample the curves $NX$ and $NK$ with 5 points per
unit length. $\{f_j, j\in J\}$ are randomly generated numbers with
mean $0$. Suppose we use $\{u^a_i, i\in I\}$ to denote the results of
our algorithm. To study the accuracy of our algorithm, we pick a set
$S \subset I$ of size 200 and estimate the error by
\[
\sqrt{
\frac{ \sum_{i\in S} |u_i - u_i^a |^2 } { \sum_{i\in S} |u_i|^2 }
}
\]
where $\{u_i, i\in S\}$ are the exact potentials computed by direct
evaluation.

Before reporting the numerical results, let us summarize the notations
that are used in the tables. $N$ is the size of the domain, $p$ is the
size of the Cartesian grid used for the equivalent charges, $P$ is the
maximum of the numbers of points in $\{x_i\}$ and $\{k_j\}$, $T_a$ is
the running time of our algorithm in seconds, $T_d$ is the estimated
running time of the direct evaluation in seconds, $T_d/T_a$ is the
speedup factor, and finally $\eps_a$ is the approximation error.

\begin{table}[h]
  \begin{center}
    \includegraphics[height=2in]{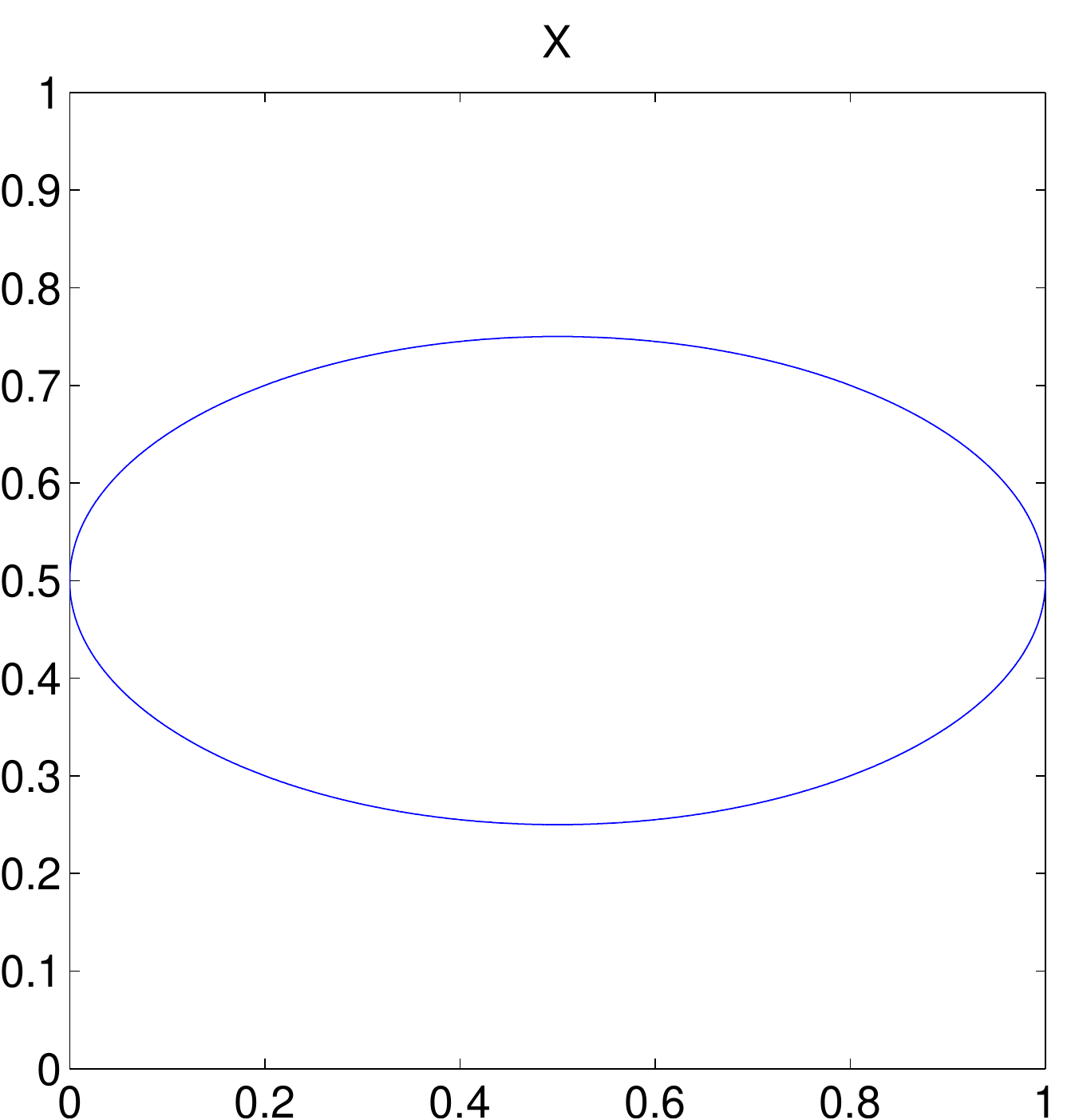}
    \includegraphics[height=2in]{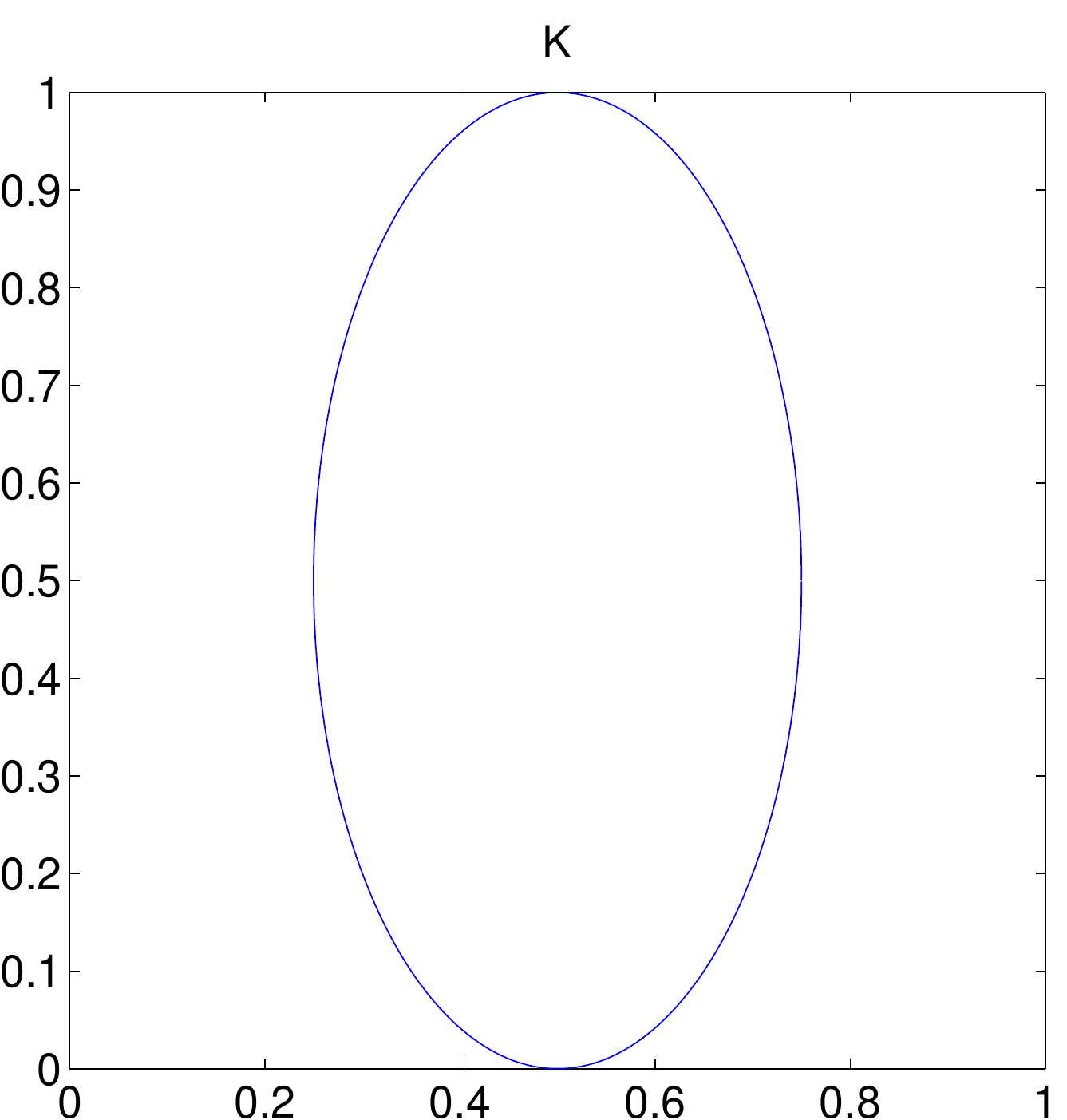}\\
    \vspace{0.1in}
    \begin{tabular}{|cccccc|}
      \hline
      $(N,p)$ & $P$ & $T_a$(sec) & $T_d$(sec) & $T_d/T_a$ & $\eps_a$\\
      \hline
      (1024,5) & 1.64e+4 & 1.23e+0 & 3.52e+1 & 2.86e+1 & 1.66e-3\\
      (2048,5) & 3.28e+4 & 2.71e+0 & 1.46e+2 & 5.38e+1 & 1.76e-3\\
      (4096,5) & 6.55e+4 & 6.05e+0 & 5.77e+2 & 9.53e+1 & 1.94e-3\\
      (8192,5) & 1.31e+5 & 1.32e+1 & 2.35e+3 & 1.78e+2 & 1.97e-3\\
      (16384,5) & 2.62e+5 & 2.89e+1 & 9.40e+3 & 3.25e+2 & 2.00e-3\\
      (32768,5) & 5.24e+5 & 6.27e+1 & 3.76e+4 & 5.99e+2 & 2.11e-3\\
      \hline
      (1024,7) & 1.64e+4 & 2.07e+0 & 3.60e+1 & 1.74e+1 & 9.31e-6\\
      (2048,7) & 3.28e+4 & 4.64e+0 & 1.44e+2 & 3.11e+1 & 9.20e-6\\
      (4096,7) & 6.55e+4 & 1.03e+1 & 5.83e+2 & 5.68e+1 & 1.03e-5\\
      (8192,7) & 1.31e+5 & 2.27e+1 & 2.35e+3 & 1.04e+2 & 1.04e-5\\
      (16384,7) & 2.62e+5 & 5.14e+1 & 9.40e+3 & 1.83e+2 & 1.07e-5\\
      (32768,7) & 5.24e+5 & 1.18e+2 & 3.76e+4 & 3.18e+2 & 1.18e-5\\
      \hline
      (1024,9) & 1.64e+4 & 3.31e+0 & 3.60e+1 & 1.09e+1 & 4.77e-8\\
      (2048,9) & 3.28e+4 & 7.23e+0 & 1.44e+2 & 1.99e+1 & 5.85e-8\\
      (4096,9) & 6.55e+4 & 1.59e+1 & 5.80e+2 & 3.65e+1 & 5.05e-8\\
      (8192,9) & 1.31e+5 & 3.57e+1 & 2.35e+3 & 6.59e+1 & 5.75e-8\\
      (16384,9) & 2.62e+5 & 7.74e+1 & 9.40e+3 & 1.21e+2 & 6.16e-8\\
      (32768,9) & 5.24e+5 & 1.87e+2 & 3.76e+4 & 2.01e+2 & 5.94e-8\\
      \hline
    \end{tabular}
  \end{center}
  \caption{2D results. Top: Both $X$ and $K$ are ellipses in unit box $[0,1]^2$.
    Bottom: Running time, speedup factor and accuracy for different combinations
    of $N$ and $p$.  $N$ is the size of the domain, $p$ is the
    size of the Cartesian grid used for the equivalent charges, $P$ is the
    maximum of the numbers of points in $\{x_i\}$ and $\{k_j\}$, $T_a$ is
    the running time of our algorithm in seconds, $T_d$ is the estimated
    running time of the direct evaluation in seconds, $T_d/T_a$ is the
    speedup factor, and finally $\eps_a$ is the approximation error.
  }
  \label{tbl:2dellp}
\end{table}

\begin{table}[h]
  \begin{center}
    \includegraphics[height=2in]{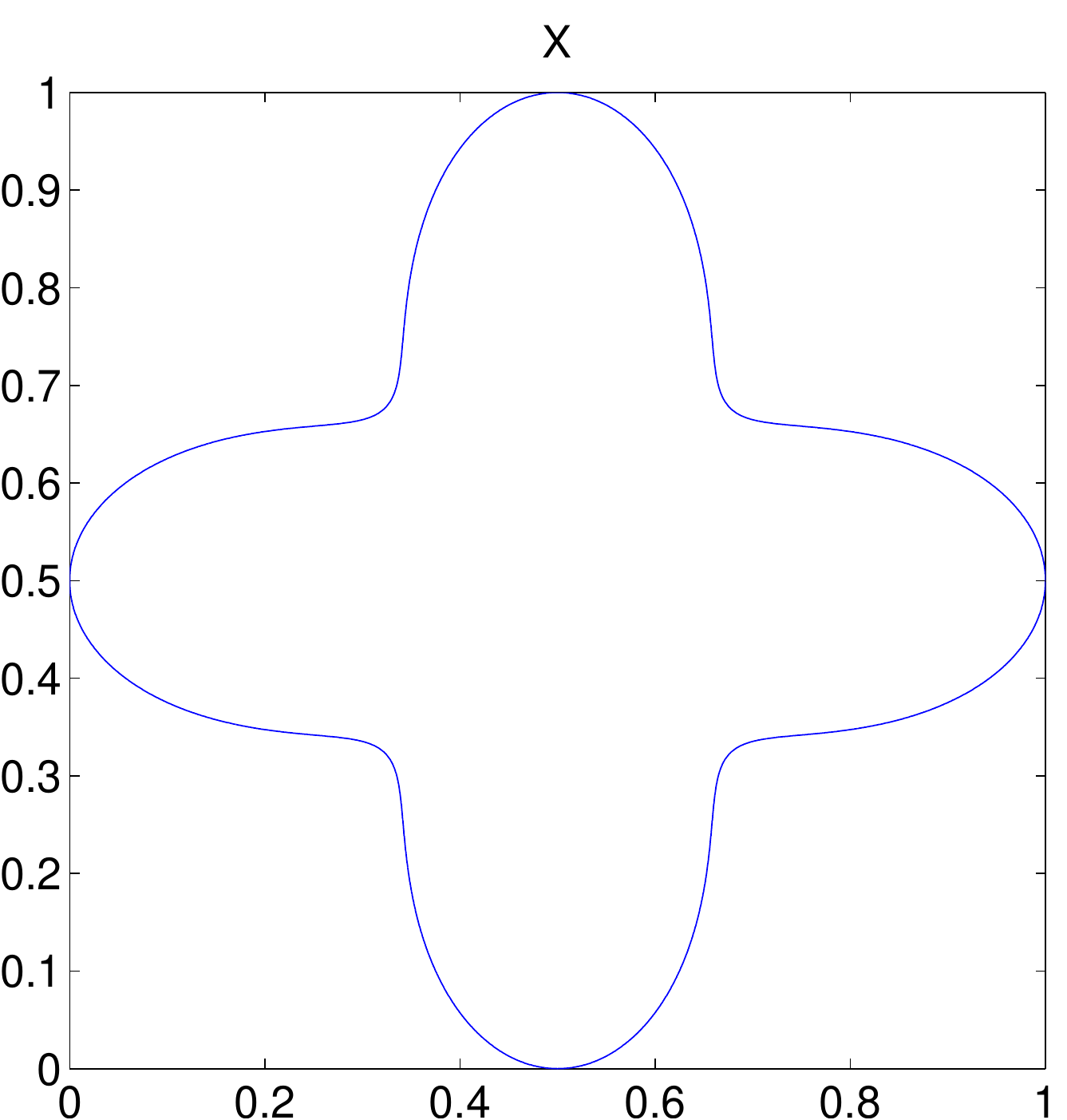}
    \includegraphics[height=2in]{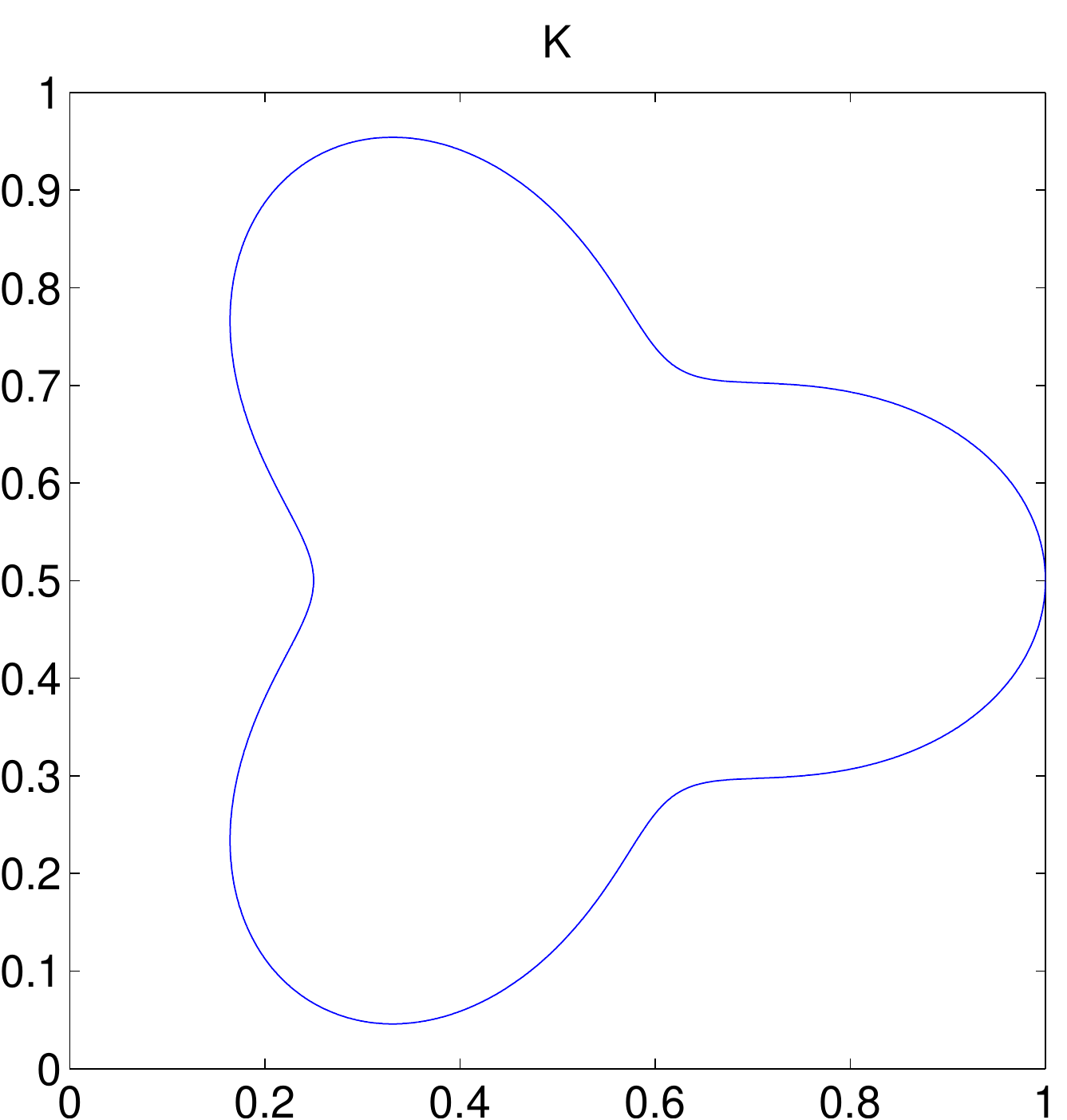}\\
    \vspace{0.1in}
    \begin{tabular}{|cccccc|}
      \hline
      $(N,p)$ & $P$ & $T_a$(sec) & $T_d$(sec) & $T_d/T_a$ & $\eps_a$\\
      \hline
      (1024,5) & 1.64e+4 & 1.93e+0 & 3.60e+1 & 1.87e+1 & 1.52e-3\\
      (2048,5) & 3.28e+4 & 4.25e+0 & 1.46e+2 & 3.43e+1 & 1.67e-3\\
      (4096,5) & 6.55e+4 & 9.49e+0 & 5.80e+2 & 6.11e+1 & 1.66e-3\\
      (8192,5) & 1.31e+5 & 2.06e+1 & 2.34e+3 & 1.14e+2 & 1.69e-3\\
      (16384,5) & 2.62e+5 & 4.59e+1 & 9.50e+3 & 2.07e+2 & 1.99e-3\\
      (32768,5) & 5.24e+5 & 9.85e+1 & 3.78e+4 & 3.84e+2 & 1.84e-3\\
      \hline
      (1024,7) & 1.64e+4 & 3.20e+0 & 3.60e+1 & 1.13e+1 & 7.81e-6\\
      (2048,7) & 3.28e+4 & 7.12e+0 & 1.44e+2 & 2.02e+1 & 8.26e-6\\
      (4096,7) & 6.55e+4 & 1.60e+1 & 5.77e+2 & 3.60e+1 & 8.79e-6\\
      (8192,7) & 1.31e+5 & 3.54e+1 & 2.34e+3 & 6.61e+1 & 9.25e-6\\
      (16384,7) & 2.62e+5 & 8.14e+1 & 9.52e+3 & 1.17e+2 & 9.07e-6\\
      (32768,7) & 5.24e+5 & 1.91e+2 & 3.76e+4 & 1.97e+2 & 1.07e-5\\
      \hline
      (1024,9) & 1.64e+4 & 5.02e+0 & 3.52e+1 & 7.02e+0 & 4.24e-8\\
      (2048,9) & 3.28e+4 & 1.12e+1 & 1.43e+2 & 1.28e+1 & 4.77e-8\\
      (4096,9) & 6.55e+4 & 2.47e+1 & 5.87e+2 & 2.37e+1 & 4.65e-8\\
      (8192,9) & 1.31e+5 & 5.60e+1 & 2.33e+3 & 4.17e+1 & 4.35e-8\\
      (16384,9) & 2.62e+5 & 1.24e+2 & 9.40e+3 & 7.60e+1 & 4.99e-8\\
      (32768,9) & 5.24e+5 & 2.84e+2 & 3.76e+4 & 1.32e+2 & 6.04e-8\\
      \hline
    \end{tabular}
  \end{center}
  \caption{
    2D results. Top: $X$ and $K$ are two smooth curves in unit box $[0,1]^2$.
    Bottom: Running time, speedup factor and accuracy for different combinations
    of $N$ and $p$.
  }
  \label{tbl:2dmult}
\end{table}

Tables \ref{tbl:2dellp} and \ref{tbl:2dmult} report the results
for two testing examples. From these tables, it is quite clear that
the complexity of our algorithm grows indeed almost linearly in terms
of the number of points, and its accuracy is stably controlled by the
value of $p$. For larger values of $N$, we obtain a substantial
speedup over the direct evaluation.

\subsection{3D case}
We apply our algorithm to the problem of computing the far field
pattern \eqref{eq:far}. In this setup, $X$ is always a sphere, while
$K$ is the boundary of the scatter. We sample the surface $NX$ and
$NK$ again with about $5 \times 5 = 25$ points per unit area.  Tables
\ref{tbl:3dplane} and \ref{tbl:3dsub} summarize the results of two
typical examples in evaluating the far field pattern of scattering
fields.

\begin{table}[h]
  \begin{center}
    \includegraphics[height=2.5in]{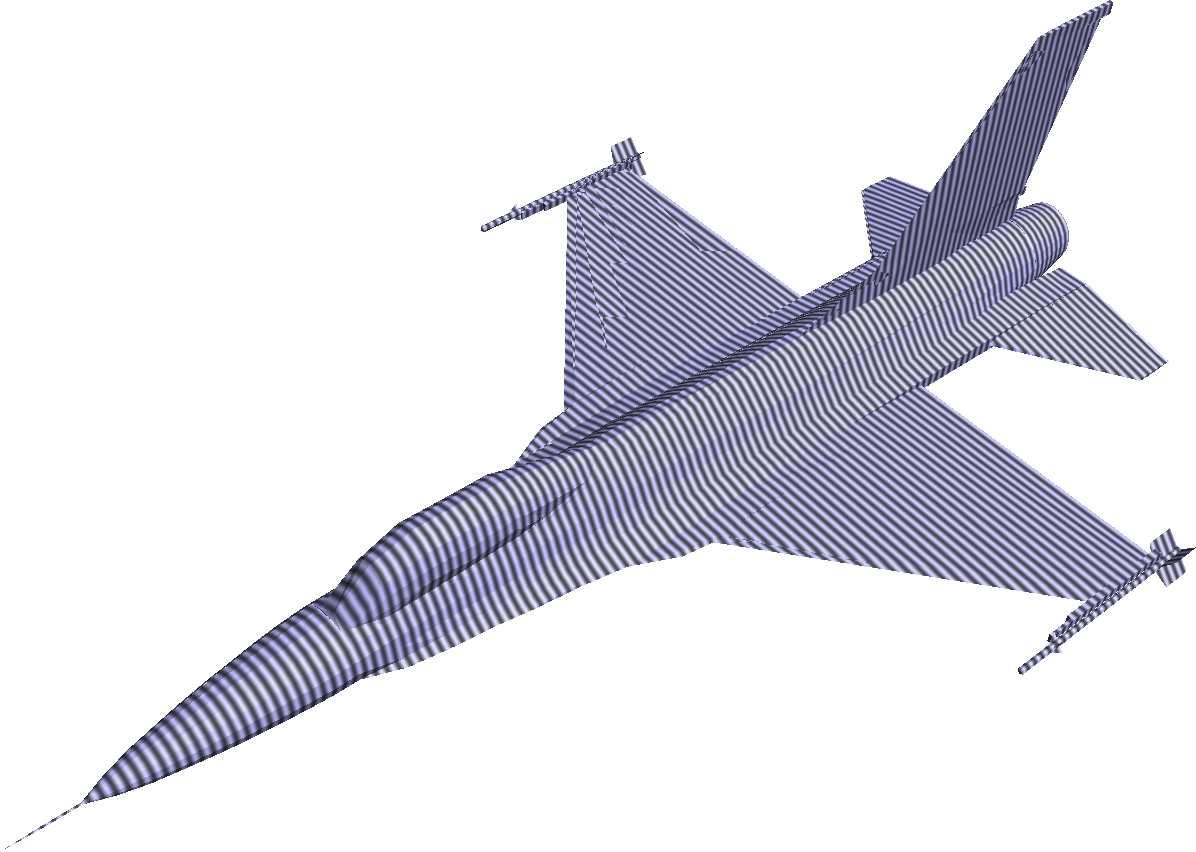}\\
    \vspace{0.1in}
    \begin{tabular}{|cccccc|}
      \hline
      $(N,p)$ & $P$ & $T_a$(sec) & $T_d$(sec) & $T_d/T_a$ & $\eps_a$\\
      \hline
      (16,5) & 2.14e+4 & 2.49e+0 & 1.50e+1 & 6.03e+0 & 1.24e-3\\
      (32,5) & 8.19e+4 & 9.78e+0 & 1.88e+2 & 1.93e+1 & 1.52e-3\\
      (64,5) & 3.22e+5 & 3.94e+1 & 2.77e+3 & 7.03e+1 & 1.44e-3\\
      (128,5) & 1.28e+6 & 1.62e+2 & 4.39e+4 & 2.71e+2 & 1.68e-3\\
      (256,5) & 5.13e+6 & 6.77e+2 & 7.00e+5 & 1.03e+3 & 1.79e-3\\
      \hline
      (16,7) & 2.14e+4 & 6.77e+0 & 1.50e+1 & 2.22e+0 & 5.80e-5\\
      (32,7) & 8.19e+4 & 2.66e+1 & 1.80e+2 & 6.79e+0 & 7.24e-5\\
      (64,7) & 3.22e+5 & 1.07e+2 & 2.74e+3 & 2.55e+1 & 7.98e-5\\
      (128,7) & 1.28e+6 & 4.40e+2 & 4.39e+4 & 9.98e+1 & 7.89e-5\\
      \hline
      (16,9) & 2.14e+4 & 1.43e+1 & 1.50e+1 & 1.05e+0 & 2.40e-7\\
      (32,9) & 8.19e+4 & 5.64e+1 & 1.88e+2 & 3.34e+0 & 3.25e-7\\
      (64,9) & 3.22e+5 & 2.28e+2 & 2.74e+3 & 1.20e+1 & 3.24e-7\\
      (128,9) & 1.28e+6 & 9.38e+2 & 4.44e+4 & 4.74e+1 & 3.33e-7\\
      \hline
    \end{tabular}
  \end{center}
  \caption{3D results. Top: the surface $K$ is the boundary of an F16 airplane model.
    Bottom: Running time, speedup factor and accuracy for different combinations
    of $N$ and $p$.  }
  \label{tbl:3dplane}
\end{table}

\begin{table}[h]
  \begin{center}
    \includegraphics[height=2.5in]{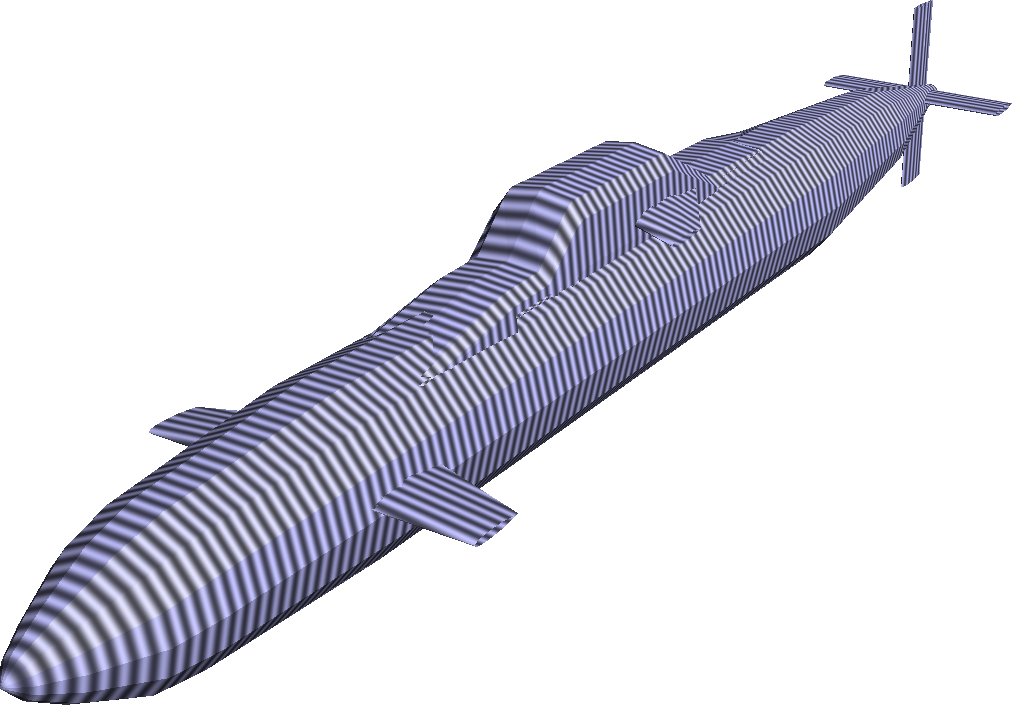}\\
    \vspace{0.1in}
    \begin{tabular}{|cccccc|}
      \hline
      $(N,p)$ & $P$ & $T_a$(sec) & $T_d$(sec) & $T_d/T_a$ & $\eps_a$\\
      \hline
      (16,5) & 2.14e+4 & 3.20e+0 & 1.07e+1 & 3.35e+0 & 1.45e-3\\
      (32,5) & 8.19e+4 & 1.25e+1 & 1.39e+2 & 1.12e+1 & 1.65e-3\\
      (64,5) & 3.22e+5 & 5.11e+1 & 2.13e+3 & 4.17e+1 & 1.79e-3\\
      (128,5) & 1.28e+6 & 2.02e+2 & 3.46e+4 & 1.71e+2 & 2.19e-3\\
      (256,5) & 5.13e+6 & 8.31e+2 & 5.54e+5 & 6.67e+2 & 1.94e-3\\
      \hline
      (16,7) & 2.14e+4 & 8.72e+0 & 1.07e+1 & 1.23e+0 & 7.06e-5\\
      (32,7) & 8.19e+4 & 3.39e+1 & 1.39e+2 & 4.11e+0 & 7.57e-5\\
      (64,7) & 3.22e+5 & 1.35e+2 & 2.13e+3 & 1.57e+1 & 9.05e-5\\
      (128,7) & 1.28e+6 & 5.48e+2 & 3.44e+4 & 6.28e+1 & 1.04e-4\\
      \hline
      (16,9) & 2.14e+4 & 1.85e+1 & 8.58e+0 & 4.64e-1 & 2.61e-7\\
      (32,9) & 8.19e+4 & 7.19e+1 & 1.39e+2 & 1.94e+0 & 3.00e-7\\
      (64,9) & 3.22e+5 & 2.87e+2 & 2.13e+3 & 7.41e+0 & 4.38e-7\\
      (128,9) & 1.28e+6 & 1.17e+3 & 3.52e+4 & 3.01e+1 & 3.46e-7\\
      \hline
    \end{tabular}
  \end{center}
  \caption{3D results. Top: the surface $K$ is the boundary of a submarine model.
    Bottom: Running time, speedup factor and accuracy for different combinations
    of $N$ and $p$.  }
   \label{tbl:3dsub}
\end{table}

\section{Conclusions and Discussions}
\label{sec:conc}

In this paper, we introduced an efficient algorithm for computing
sparse Fourier transforms located on curves and surfaces. Our
algorithm, which is an extension of the butterfly algorithm, is
accurate and has provably $O(N \log N)$ complexity. The success of the
algorithm is based on an low rank property concerning the interaction
between spatial and frequency regions that follow a certain
geometrical condition. We use equivalent sources supported on
Cartesian grids as the low rank representation and exploit the
tensor-product property of the Fourier transform to achieve maximum
efficiency.  Furthermore, our algorithm requires only linear storage
space.

The problem considered in this paper is only one example of many
computational issues regarding highly oscillatory behaviors. Some
other examples include the computation of Fourier integer operators
\cite{stein-1993-ha}, scattering fields for high frequency waves
\cite{colton-1983-iemst}, and Kirchhoff migrations
\cite{beylkin-1984-ipagrt}. In the past two decades, many algorithms
have been developed to address these challenging computational tasks
efficiently. Some examples include
\cite{bruno-2001-afho,cheng-2006-awfmm,averbuch-2000-ecoi,candes-2007-fcfio,engquist-2007-fdmaok}.
It would be interesting to see whether the ideas behind these
approaches can be used to study the problem addressed in this paper,
and vice versa.

{\bf Acknowledgments.} The research presented in this paper was
supported by an Alfred P. Sloan Fellowship and a startup grant from
the University of Texas at Austin.

\bibliographystyle{abbrv}
\bibliography{ref}

\end{document}